\documentclass[11pt]{amsart}
\usepackage{geometry}                
\usepackage{graphicx}
\usepackage{float}
\usepackage{amssymb}
\usepackage{epstopdf}
\usepackage{color}
\usepackage{tikz-cd}
\usepackage{capt-of}

\newtheorem{thm}{Theorem}

\newtheorem{prop}[thm]{Proposition}
\newtheorem{lem}[thm]{Lemma}

\usepackage{cite}
\usepackage[backref]{hyperref}

 \usepackage[usenames,dvipsnames]{pstricks}
 \usepackage{pstricks-add}
 \usepackage{epsfig}
 \usepackage{pst-grad} 
 \usepackage{pst-plot} 
 \usepackage[space]{grffile} 
 \usepackage{etoolbox} 
 \makeatletter 
 \patchcmd\Gread@eps{\@inputcheck#1 }{\@inputcheck"#1"\relax}{}{}
 \makeatother

\title{A study on Diophantine equations via cluster theory}
\author{Leizhen Bao \; and \; Fang Li}

\address{Leizhen Bao
\newline School of Mathematical Sciences, 
Zhejiang University,  
Hangzhou, Zhejiang 310058, 
China P.R.}
\email{mathbao@zju.edu.cn}

   \address{Fang Li
\newline School of Mathematical Sciences, 
Zhejiang University,  
Hangzhou, Zhejiang 310058, 
China P.R.}
\email{fangli@zju.edu.cn}

       
\date{version of \today}                                   

\begin{document}

\maketitle

\begin{abstract}

In this paper, we mainly answer a Lampe's question\cite{lampe} about the solutions of a Diophantine equation, that is, we give a criterion to determine which solutions of the Diophantine equation are in the orbit of the initial solution $(\epsilon,\epsilon,\epsilon,
\epsilon,\epsilon)$ under the actions of the group $G$ which is defined by mutations of a cluster algebra. In order to do this, using a rational map $\varphi$, we transform the Diophantine equation to a related equation whose all positive integral solutions form the orbit of an initial solutions $\varphi(\epsilon,\epsilon,\epsilon,
\epsilon,\epsilon) = (3,4,4)$ under the actions of the group $\widetilde{G}$, and the set $S(3,4,4)$ is shown to be the orbit of $(\epsilon,\epsilon,\epsilon,
\epsilon,\epsilon)$ under the actions of a subgroup of $G$. Then the criterion is proved as the main conclusion.
\end{abstract}

\tableofcontents

\section{Introduction}

\subsection{Motivation}

In 1880, when Markov\cite{Markov} study the following Diophantine equation (called \textbf{Markov equation}) 
	$$T'(X,Y,Z) := {X^2 + Y^2 + Z^2}-3{XYZ} = 0,$$
he discovered that the set of whole positive integer solutions of this equation is the orbit of the initial solution $(1,1,1)$ under the group $G' := \langle\mu'_1, \mu'_2, \mu'_3\rangle$, where 
\begin{equation}\label{eq: makrov mu}
	\mu'_1(x,y,z) = (3yz-x,y,z), \mu'_2(x,y,z) = (x, 3xz-y, z), \mu'_3(x,y,z) = (x,y, 3xy - z).
\end{equation}
That is, we have\footnote{For convenience,   we write the orbit of $(x_1, ..., x_n)$ under the action of the group $G$ as $G(x_1, ..., x_n):= \{g(x_1, ..., x_n) \mid \forall g \in G \} $.} following relation
\begin{equation*}
		G'(1,1,1) = \{ (x,y,z) \in \mathbb{N} ^3 \mid T'(x,y,z) = 0 \}.
\end{equation*}

	In 2002, for studying problems in representation theory, Fomin and Zelevinsky\cite{FZ1} abstracted out an algebraic structure, named \textit{cluster algebra}. Using an integer matrix $B$ and variables $\mathbf{x}$ via transformations called \textit{mutations}, we can define a cluster algebra $\mathcal{A}(\mathbf{x}, B)$. These definitions will be recalled in section \ref{sec:ca}. 
	
	 In 2012, Peng and Zhang\cite{PengZhang} revealed the connection between the Markov equation and a cluster algebra. Denoting a matrix 
\begin{equation*} B' := \left(\begin{array}{rrr}
  0 & 2 & -2 \\
  -2 & 0 & 2 \\
  2 & -2 & 0
\end{array} \right),
\end{equation*}
they found  the mutations of the cluster algebra $\mathcal{A}((x_1,x_2,x_3), B')$ are just the actions $\mu'_1, \mu'_2, \mu'_3$ defined in \eqref{eq: makrov mu}.\\

	On the other hand, in 2016 Lampe\cite{lampe} found that if we let a matrix
    \begin{equation*} B'' := \left(\begin{array}{rrr}
  0 & 1 & -1 \\
 -4 & 0 & 2 \\
  4 & -2 & 0
\end{array} \right),
\end{equation*}
then there arose a variant of the Markov equation
	$$T''(X,Y,Z) := {X^2 + Y^4 + Z^4 + 2XY^2 + 2XZ^2}-7{XY^2Z^2} = 0.$$
	The set of all positive integer solutions of this equation is the orbit of the initial solution $(1,1,1)$ under the  group $G'' := \langle\mu''_1, \mu''_2, \mu''_3\rangle$, where $\mu''_i$'s are mutation of cluster algebra $\mathcal{A}((x_1,x_2,x_3), B'')$, defined by 
	$$\mu''_1(x,y,z) = (\frac{y^4+z^4}{x}, y, z), \mu''_2(x,y,z) = (x, \frac{x+z^2}{y}, z), \mu''_3(x,y,z) = (x, y,\frac{x+y^2}{z}).$$ Which means, we have
\begin{equation*}
	G''(1,1,1) = \{ (x,y,z) \in \mathbb{N} ^3 \mid T''(x,y,z) = 0 \}.\\
\end{equation*}

And in the same article, Lampe introduced another cluster algebra $\mathcal{A}((x_1,x_2,x_3,x_4,x_5), B)$, where
\begin{equation*}
		{B} := \left(\begin{array}{rrrr}
		0 & -2 & 1 & 1\\
		2 & 0 & -1 & -1\\
		-1 & 1 & 0  & 1\\
		-1 & 1 & -1 & 0\\
		0 & 0 & 1 & -1
\end{array}\right).
\end{equation*}
Using the mutations of this cluster algebra, he defined a group $G:=\langle\alpha, \beta\rangle$, where {$\alpha, \beta$ are maps from $\mathbb{Q}^5_+$ to $\mathbb{Q}^5_+$ defined as
\begin{equation*}
	\alpha(a,b,c,d,e) :=  ( b,\frac{b^2+cd}{a}, c,d,e),\quad \beta(a,b,c,d,e) := ( b,c, \frac{ac+be}{d}, a, e).
\end{equation*}
Clearly, they are invertible.}

Then Lampe discovered that for any action $g \in G$, the vector $g(1,1,1,1,1)$ is a positive integer solution of the Diophantine equation
\footnote{  The equation \eqref{T=9} is different with the original equation defined by Lampe \cite{lampe} via subtracting $9$ in the two-sides for the convenience of calculation. }
  $T(a,b,c,d,1) = 0$, where
\begin{equation}\label{T=9}
		T(a,b,c,d,e) := \frac{ab(c^2+d^2+e^2) + (a^2+b^2+cd)(c+d)e}{abcd}-9.
	\end{equation}
What is more, for any action $g \in G$ and any constant $\epsilon \in \mathbb{N}$, the vector $g(\epsilon, \epsilon, \epsilon, \epsilon, \epsilon)$ is a positive integer solution of the Diophantine equation $T(a,b,c,d,\epsilon) = 0$. In the opposite direction, however, by computer calculations suggest that not all solutions can be obtained by this way. Which means, for any $\epsilon \in \mathbb{N}$, we have
$$G(\epsilon, \epsilon, \epsilon, \epsilon, \epsilon) \subsetneq \{(a,b,c,d, \epsilon) \in \mathbb{N}^5 \mid T(a,b,c,d, \epsilon) = 0 \}.$$

	Therefore, Lampe raised following questions:\vspace{2mm}
\\
	\textbf{Question (i).} Which solutions $(a,b,c,d,1) \in \mathbb{N}^5$ of the Diophantine equation $T(a,b,c,d,1) = 0$ can be obtained from the initial solution $(1,1,1,1,1)$ by a sequence of actions of $G$? \vspace{2mm}
	\\
	\textbf{Question (ii).} Fix $\epsilon\in\mathbb{N}$. Which solutions $(a,b,c,d,\epsilon) \in \mathbb{N}^5$ of the Diophantine equation $T(a,b,c,d,\epsilon) = 0$ can be obtained from the initial solution $(\epsilon,\epsilon,\epsilon,\epsilon,\epsilon)$ by a sequence of actions of $G$?  \vspace{2mm}

	In this article, we will answer above questions.
	We first introduce the cluster algebra $\mathcal{A}((a,b,c,d, \epsilon), B)$ and the Diophantine equation $T(a,b,c,d, \epsilon) = 0$, and study some basic  propositions of them in section \ref{sec:lampe}. 
	Then, we will directly answer the question (ii). The question (i) is just a special case of (ii). We solve it by three steps.
	
	STEP  1.
 
 In section \ref{sec:hone}, we introduce a rational map $\varphi: \mathbb{Q}_+^5 \to \mathbb{Q}_+^3$, and find that for any solution $P$ of the equation $T(P)=0$, the triple $\varphi(P)$ can be obtained from the triple $(3,4,4)$ by a sequence of actions of the group $\widetilde{G}$ which is 3-dimension version of the group $G$. 
	
	STEP  2.
 
 In section \ref{sec:pell}, we prove that every vector $P$, which satisfies $T(P)=0$ and $\varphi(P)=(3,4,4)$, can be obtained from the initial solution $(\epsilon,\epsilon,\epsilon,\epsilon,\epsilon)$ by a sequence of actions of the group $G$.
	
	STEP  3. 
 
 In section \ref{sec:bao}, we prove that the solution $P = (a,b,c,d,\epsilon)$, which satisfying $\varphi(P) \in \mathbb{N}^3$ and $P \equiv 0 (\textrm{mod~} \epsilon)$, can be obtained from the initial solution $(\epsilon,\epsilon,\epsilon,\epsilon,\epsilon)$ by a sequence of actions of the group $G$.
	
	Hence the answer of the \textbf{Question (ii)} is \text{Theorem \ref{MainThm}}, which is the main result in this paper, as follows.\\
    {\bf Theorem} \ref{MainThm}.
      {\em  Fix the constant $\epsilon \in \mathbb{N}$. The positive integer solutions $(a,b,c,d,\epsilon) \in \mathbb{N}^5$ of the Diophantine equation $T(a,b,c,d,\epsilon)=0$ defined in \eqref{T=9} can be obtained from the initial solution $(\epsilon,\epsilon,\epsilon,\epsilon,\epsilon)$ by a sequence of actions of the group $G:= \langle\alpha,\beta\rangle$ defined in \eqref{eq: alpha beta}, if and only if, $\varphi(a,b,c,d,\epsilon) \in \mathbb{N}^3 \textrm{~and~} (a,b,c,d,\epsilon) \equiv 0 (\textrm{mod~}\epsilon)$, where the map $\varphi$  defined in \eqref{varphi}. That is, we have 
	$$G(\epsilon,\epsilon,\epsilon,\epsilon,\epsilon) = \{ P=(a,b,c,d,\epsilon) \in \mathbb{N}^5 \mid T(P) = 0, \varphi(P) \in \mathbb{N}^3, P \equiv 0
	(\textrm{mod~}\epsilon) \}.$$  }

\textbf{Remark.} For convenience, in this paper we write $(x_1, \dots, x_n) \equiv r (\textrm{mod~} p)$ if $x_i \equiv r (\textrm{mod~} p)$ for $i = 1, 2, \dots , n$. In our context $\epsilon \in \mathbb{N}$ is a constant, so we fix it and without any mention in following paragraph.

\subsection{What are cluster algebras?}\label{sec:ca}

We first recall some definitions of cluster algebra. 
We fix $m, n \in \mathbb{N}$ satisfying $m \geq n$. For brevity, we write ${I} = \{1, 2, \dots, m\}$ and ${J} = \{1, 2, \dots, n\}$. Given an integer matrix ${B} = (b_{ij})_{i\in {I}, j \in {J}}$, we called it \textit{extended exchange matrix}\footnote{When $m=n$, we call it \textit{exchange matrix}.}. The \textit{mutation} of matrix ${B}$ in the direction $k \in {J}$ is a new matrix $\mu_k({B}) =  {B}'  = (b_{ij}')_{i\in {I}, j \in {J}}$, 
where
\begin{equation*}
	b_{ij}' = \begin{cases} -b_{ij},  & \textrm{if~} i=k \textrm{~or~} j =k;\\
	b_{ij} + \frac{1}{2}(b_{ik}|b_{kj}| + |b_{ik}|b_{kj}), & \textrm{otherwise}.
\end{cases}
\end{equation*}

An \textit{extended cluster} is a  $m$-tuple ${\mathbf{x}} = (x_1, \dots, x_m)$, where $x_1, \dots, x_m$ are algebraically independent variables over $\mathbb{Q}$. 
	The $n$-tuple ${\mathbf{x}_0} = (x_1, \dots, x_n)$ is called \textit{cluster}, and its elements $x_1, \dots, x_n$ are called \textit{cluster variables}. The other elements $x_{n+1}, \dots, x_m$ are called \textit{frozen variables}. 
	The \textit{mutation} of extended cluster ${\mathbf{x}}$  in the direction $k \in {J}$ is a new extended cluster $\mu_k({\mathbf{x}}) = {\mathbf{x}}'$, where
	${\mathbf{x}}' = (x_1', \dots, x_m')$ is given by $x_i' = x_i$ for $i \ne k$, whereas $x_k'$ is determined by the \textit{exchange relation}
\begin{equation*}
	x_kx_k' = \prod_{i \in {I}, b_{ik}>0}x_i^{b_{ik}} + \prod_{i \in {I}, b_{ik}<0}x_i^{-b_{ik}}.
\end{equation*}

	A \textit{seed} is a pair $({\mathbf{x}}, {B})$, where ${\mathbf{x}}$ is a extended cluster and ${B}$ is an integer $m \times n$ matrix whose \textit{principal part} $B_0=(b_{ij})_{i,j \in {J}}$ is \textit{skew-symmetrizable} matrix, i.e. there exists a diagonal matrix $D=diag(d_1, \dots, d_n)$ with $d_i > 0$ such that the matrix $DB_0$ is skew-symmetric. 
	The \textit{mutation} of seed $({\mathbf{x}}, {B})$  in the direction $k \in {J}$ is a new seed $\mu_k({\mathbf{x}}, {B}) =  ({\mathbf{x}}', {B}')$, where  ${B}' = \mu_k({B})$ and ${\mathbf{x}}' = \mu_k({\mathbf{x}})$.
	It is easy to check that mutation is involutive, i.e. we have $\mu_k^2({\mathbf{x}}, {B}) = ({\mathbf{x}}, {B})$. 
	
	The \textit{permutation} of seed $({\mathbf{x}}, {B})$ is a permutation $\sigma$ of ${I}$ which preserves ${I} \setminus {J}$ setwise. It acts on a seed as $\sigma({\mathbf{x}}, {B}) = ({\mathbf{x}}', {B}')$, where
	$$b_{ij}' = b_{\sigma^{-1}{(i)}\sigma^{-1}{(j)}}, \quad x_i' = x_{\sigma^{-1}{(i)}}.$$
In this case, we called the two seeds are \textit{isomorphic}, and writes as $({\mathbf{x}}, {B}) \simeq (\tilde{\mathbf{x}'}, {B}')$.
	Two seeds $({\mathbf{x}}, {B})$ and $({\mathbf{x}}', {B}')$ are called \textit{mutation equivalent} if there is a sequence of mutations taking one to the other up to isomorphic, i.e. there exists a sequence $(k_1, \dots, k_r) \in {J}^r$ of length $r \geq 0$ such that $\mu_{k_r} \cdots \mu_{k_1}({\mathbf{x}}, {B}) \simeq ({\mathbf{x}}', {B}')$. In this case, we write $({\mathbf{x}}, {B}) \sim ({\mathbf{x}}', {B}')$.

	Now we can define cluster algebra. Let the pair $({\mathbf{x}}, {B})$ be a seed. The field $\mathcal{F}=\mathbb{Q}(x_1, \cdots, x_m)$ is called \textit{ambient field} and the ring $\mathcal{R}=\mathbb{Q}[x_{n+1}, \cdots, x_m]$ is called \textit{ground ring}. Let
\begin{equation*}
	\mathcal{X}({\mathbf{x}}, {B}) := \bigcup_{({\mathbf{x}}, {B}) \sim ({\mathbf{x}}', {B}')} \{x_1', \dots, x_n'\}
\end{equation*}
be the set of all cluster variables appearing in its seeds.  The \textit{cluster algebra} with initial seed $({\mathbf{x}}, {B})$, denoted $\mathcal{A}({\mathbf{x}}, {B})$, is the $\mathcal{R}-$subalgebra of ambient field $\mathcal{F}$ generated by $\mathcal{X}({\mathbf{x}}, {B})$. In other words, $\mathcal{A}({\mathbf{x}}, {B}) = R[\mathcal{X}({\mathbf{x}}, {B})]$.

 	The main theorem about cluster algebras is the Laurent Phenomena\cite{FZ1}.	
\begin{thm}[Fomin, Zelevinsky\cite{FZ1}]\label{LP} 
	In  a cluster algebra $\mathcal{A}({\mathbf{x}}, {B})$, each cluster variable can be expressed as a Laurent polynomial with positive integer coefficients in the elements of any extended cluster ${\mathbf{x}}$.
\end{thm} 

    Now we explain  the quiver corresponding to a cluster algebra.
    
	When the principal part $B_0$ of ${B}$ is skew-symmetric integer matrix, we can visualize $B_0$ as a quiver.
	A \textit{quiver} is a finite directed graph, it has vertices and arrows between the vertices. 
	We can visualize the seed $({\mathbf{x}}, {B})$ as a quiver $Q({\mathbf{x}}, {B})$. 
	The quiver $Q({\mathbf{x}}, {B})$ has $m$ vertices indexed by $\{1, \dots, m\}$ and  attached with cluster variables or frozen variables, i.e. each vertex $j$ is attached with variable $x_j$.
	If $b_{ij} > 0$, there are $b_{ij}$ arrows from $i$ to $j$. Thus if $b_{ij} < 0$, there are $-b_{ij} = b_{ji}$ arrows from $j$ to $i$. 

	The mutation of $B$  in the direction $k \in {J}$ can be corresponded to the 
 \textit{mutation of a quiver} in the direction $k \in {J}$, which defined as follows:
	
	(1)  If there exist $a$ arrows $i \to j$ and $b$ arrows $k \to j$, we add $ab$ arrows $i \to j$.
	
	(2) Cancel a maximal set of 2-cycles from those created in (1).
	
	(3) Reverse all arrows incident with $k$.

\section{Some basic propositions of the equation $T=0$ and the group $G$}\label{sec:lampe}

	Now we introduce the Diophantine equation $T(a,b,c,d,e) = 0$ and the group $G$ of Lampe \cite{lampe}. We begin with a seed $({\mathbf{x}}, {B})$, where ${\mathbf{x}} = (a, b, c, d,  e)$ and
\begin{equation}\label{mat: B}
	{B} = \left(\begin{array}{rrrr}
		0 & -2 & 1 & 1\\
		2 & 0 & -1 & -1\\
		-1 & 1 & 0  & 1\\
		-1 & 1 & -1 & 0\\
		0 & 0 & 1 & -1
\end{array}\right).
\end{equation}

\begin{figure}[h]
\centering
	\includegraphics{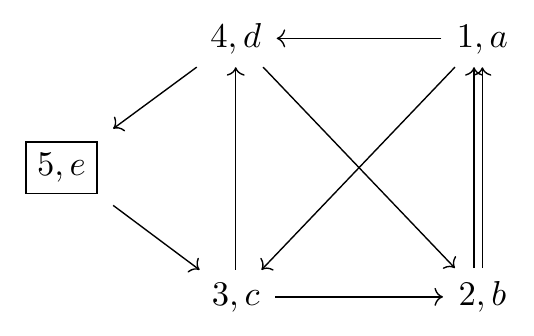}
\captionof{figure}{The quiver $Q({\mathbf{x}}, \tilde{B})$.} \label{figB}
\end{figure}

    Then the corresponding quiver is in \text{Figure \ref{figB}},

	Note that every vertex have two elements, the left one is the index of the vertex, the right one is the corresponding  cluster variable (or frozen variable). When the variable is frozen,  we box the vertex,  for example, the vertex 5 in \text{Figure \ref{figB}}.

The four mutations of ${\mathbf{x}}$ as follow,
$$\mu_1(a,b,c,d,e) = (\frac{b^2+cd}{a}, b,c,d,e), \quad \mu_2(a,b,c,d,e) = (a, \frac{a^2+cd}{b},c,d,e),$$

$$\mu_3(a,b,c,d,e) = (a,b,\frac{ae+bd}{c},d,e), \quad \mu_4(a,b,c,d,e) = (a,b, c, \frac{ac+be}{d}, e).$$

\begin{center}
\begin{figure}[h]
\centering
	\includegraphics[scale=0.85]{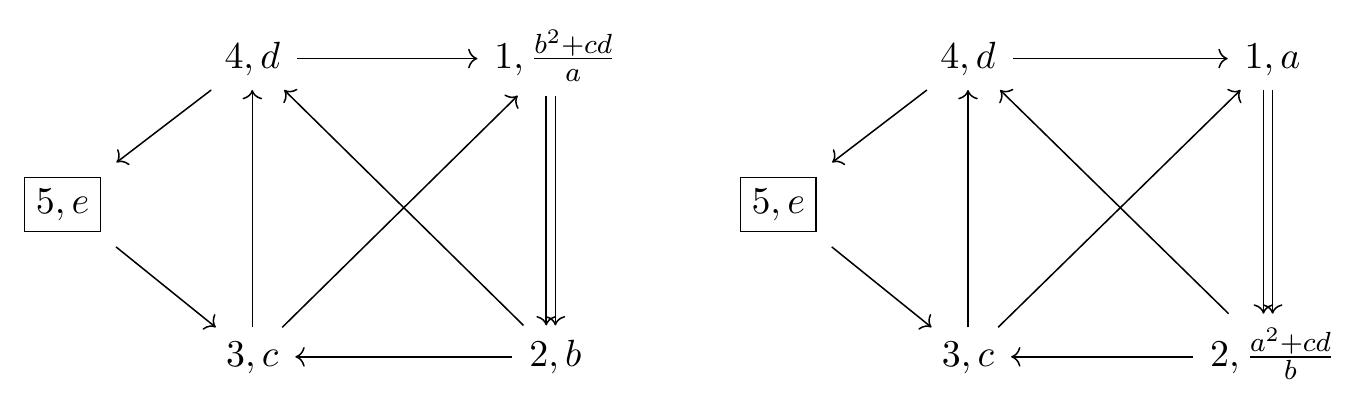}
\captionof{figure}{The quivers  $Q(\mu_1({\mathbf{x}}, \tilde{B}))$ and $Q(\mu_2({\mathbf{x}}, \tilde{B}))$.}\label{mu1mu2}
\end{figure}
\end{center}

\begin{center}
\begin{figure}[h]
\centering
	\includegraphics[scale=0.85]{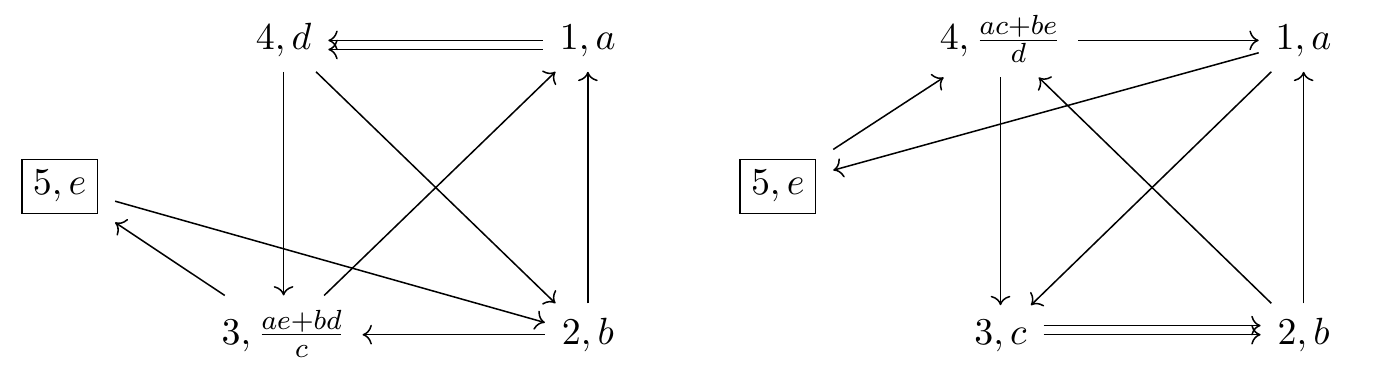}
\captionof{figure}{The quivers  $Q(\mu_3({\mathbf{x}}, \tilde{B}))$ and $Q(\mu_4({\mathbf{x}}, \tilde{B}))$.}\label{mu3mu4}
\end{figure}
\end{center}

Corresponding quivers are showed as follow Figure \ref{mu1mu2} and Figure \ref{mu3mu4},

Note that the structure of quiver does not change after any mutation, hence we have the following proposition.

\begin{prop}[{Lampe \cite{lampe}}]
	The matrix $B$ defined in \eqref{mat: B} is an invariant under the mutations of the cluster algerbra $\mathcal{A}((a,b,c,d,e), B)$ up to permutation. That is, we have 
	$${B} = \sigma_{(12)}\mu_1({B}) = \sigma_{(12)}\mu_2({B}) = \sigma_{(1234)}\mu_3({B}) = \sigma_{(4321)}\mu_4({B}),$$
where $\sigma_{(12)}, \sigma_{(1234)}, \sigma_{(4321)} \in S_5$.
\end{prop}


	Inspired by this proposition, we extend the mutation by setting
\begin{equation}
	\alpha := \sigma_{(12)}\mu_1, \quad  \beta := \sigma_{(4321)}\mu_4.
\end{equation}
	Then we have
\begin{equation}\label{eq: alpha beta}
	\alpha(a,b,c,d,e) :=  ( b,\frac{b^2+cd}{a}, c,d,e),\quad \beta(a,b,c,d,e) := (b,c,\frac{ac+be}{d},a,e).
\end{equation}
	We can check
\begin{equation*}
	\alpha^{-1} = \sigma_{(12)}\mu_2, \quad \beta^{-1} = \sigma_{(1234)}\mu_3,
\end{equation*}
	that is
\begin{equation*}
	\alpha^{-1}(a,b,c,d,e)  = (\frac{a^2+cd}{b},a,c,d,e), \quad 
	\beta^{-1}(a,b,c,d,e) = (d, a, b, \frac{ae+bd}{c}, e).
\end{equation*}


	Therefore, we define a group $G := \langle\alpha, \beta\rangle $. By  above discussion, we have
	 $$g({\mathbf{x}}, {B}) = ({\mathbf{x}}', {B}), \textrm{ for any } g \in G.$$
	Which means we can only use the group $G$ and the initial extended cluster ${\mathbf{x}}$ to compute new extended clusters without changing ${B}$.
	 
	Focusing on the group $G$, Lampe \cite{lampe} found an invariant $T$, which gives the equation \eqref{T=9}:
	\begin{equation*} 
		T(a,b,c,d,e) = \frac{ab(c^2+d^2+e^2) + (a^2+b^2+cd)(c+d)e}{abcd}-9.
	\end{equation*}

Lampe proved the following proposition.
\begin{prop}[Lampe\cite{lampe}]\label{Tg=T} 
	The map $T$ is an invariant under the actions in the group $G$.
	That is, we have 
	$$T(g(a,b,c,d,e)) = T(a,b,c,d,e),$$ 
	for all action $g \in G$ and all vector $(a,b,c,d,e) \in \mathbb{Q}_+^5$.\\
\end{prop}

    Considering the homogeneous of actions of $G$, we obtain the following proposition.
\begin{prop} \label{mod}
	The components of vectors in the orbit $G(\epsilon, \epsilon, \epsilon, \epsilon, \epsilon)$ are positive integer with multiple of $\epsilon$. That is, 
	$$G(\epsilon, \epsilon, \epsilon, \epsilon, \epsilon) \subset \{P=(a,b,c,d,\epsilon) \in \mathbb{N}^5  \mid P \equiv 0(\textrm{mod~}\epsilon)\}.$$
\end{prop}
\begin{proof}	
	For any action $g \in G$. Let a vector $P_1 = g(1,1,1,1,1)$. By \text{Theorem \ref{LP}} we know $P_1  \in \mathbb{Z}^5$. It is easy to prove that the  components of $P_1$ is positive, so  we know $P_1 \in \mathbb{N}^5$.
	
	Now we define a map $\eta_\epsilon : \mathbb{R}_+^5 \to \mathbb{R}_+^5$ by $\eta_\epsilon(a,b,c,d,e) = (a\epsilon,b\epsilon,c\epsilon,d\epsilon,e\epsilon)$. Since the action $g$ of the group $G$ is homogeneous, we have 
	    $g\eta_\epsilon = \eta_\epsilon g$.
	
	Hence by $g(\epsilon, \epsilon,\epsilon,\epsilon,\epsilon) = g(\eta_\epsilon(1,1,1,1,1)) = \eta_\epsilon(g(1,1,1,1,1)) = \eta_\epsilon(P_1)$,
	and we know that  $g(\epsilon, \epsilon,\epsilon,\epsilon,\epsilon) \in \mathbb{N}^5$ and $g(\epsilon, \epsilon,\epsilon,\epsilon,\epsilon) \equiv 0 (\textrm{mod~} \epsilon)$.
\end{proof}


\section{Solutions of the equation $\widetilde{T}=0$ under the actions of the group $\widetilde{G}$}\label{sec:hone}


In this section, we first introduce a rational map $\varphi$ which can simplify the Diophantine equation $T=0$ \eqref{T=9} to the Diophantine equation $\widetilde{T}=0$, then we will prove that the set of the positive integral solutions of this simpler equation is a orbit of actions of the group $\widetilde{G}$. We will define these notions later.
    
    We first define the rational map $\varphi: \mathbb{Q}_+^5 \to \mathbb{Q}_+^3$ by
\begin{equation} \label{varphi}
	\varphi(a, b, c, d, {{e}}) = (\mathcal{C}_0(a, b, c, d, {{e}}), \mathcal{C}_1(a, b, c, d, {{e}}), \mathcal{C}_2(a, b, c, d, {{e}})),
\end{equation}
where 
\begin{eqnarray*} 
	\mathcal{C}_0(a, b, c, d, {{e}}) &:=& \frac{a^2 + b^2 + cd}{ab},\\
	\mathcal{C}_1(a, b, c, d, {{e}}) &:=& \frac{c^2d + a^2c + b^2d+ ab{{e}}}{bcd},\\
	\mathcal{C}_2(a, b, c, d, {{e}}) &:=& \frac{cd^2 + a^2c + b^2d+ ab{{e}}}{acd}.
\end{eqnarray*}
\textbf{Remark.} The definitions of $\mathcal{C}_i$ and the following lemma are inspired by Hone's conclusion (Theorem 1 in \cite{Hone}) where he discussed a recurrence (which in our context is the action $\beta$) and corresponding conserves quantities.\\

{
    
	The map $\varphi$ can reduce the Diophantine equation $T=0$ \eqref{T=9} to a related Diophantine equation 
	\begin{equation}\label{markov+2}
	    \widetilde{T}(X,Y,Z) := XYZ - X^2 - Y^2 - Z^2 - 7.
	\end{equation}
\begin{lem}\label{TM2}
	For any vector $P = (a, b, c, d, {e}) \in \mathbb{Q}_+^5$, we have
\begin{equation}\label{T=Mvarphi+2}
	T(P) = \widetilde{T}(\varphi(P)).
\end{equation}
\end{lem}

\begin{proof} By tedious computations, we have
	\begin{eqnarray*}
	T(P) &=& \mathcal{C}_0(P)\mathcal{C}_1(P)\mathcal{C}_2(P)-\mathcal{C}_0(P)^2-\mathcal{C}_1(P)^2-\mathcal{C}_2(P)^2  -7\\
	 	&=& \widetilde{T}(\mathcal{C}_0(P),\mathcal{C}_1(P),\mathcal{C}_2(P))\\
		&=& \widetilde{T}(\varphi(P)).
	\end{eqnarray*}
\end{proof}

From the above observation, we define two mappings $\tilde{\alpha}, \tilde{\beta}$ from $\mathbb{Q}^3$ to $\mathbb{Q}^3$ for $(x,y,z) \in \mathbb{Q}^3$,
    \begin{equation*} \label{def: tab}
	    \tilde{\alpha}(x, y, z) := (x, z, xz-y), \textrm{\quad \quad} \tilde{\beta}(x, y, z) := (y, z ,x).
    \end{equation*}
Obviously, they are invertible, with inverses $\tilde{\alpha}^{-1}$ and $\tilde{\beta}^{-1}$ by
    $$\tilde{\alpha}^{-1}(x,y,z) =  (x,xy-z,y),\quad\quad \tilde{\beta}^{-1}(x,y,z) = (z, x, y).$$
So we can define $\widetilde{G}$ to be the group generated by $\tilde\alpha$ and $\tilde\beta$. That is, $\widetilde{G} := \langle\tilde{\alpha}, \tilde{\beta}\rangle$.

    The map $\varphi$ can commute with the actions of the group $G$.
 
\begin{lem}\label{varphig} Fix a vector $P \in \mathbb{Q}_+^5$. Let $g = \alpha^{n_1}\beta^{m_1} \cdots \alpha^{n_t}\beta^{m_t}$ and $\tilde{g} = \tilde{\alpha}^{n_1}\tilde{\beta}^{m_1} \cdots \tilde{\alpha}^{n_t}\tilde{\beta}^{m_t}$ for some $n_i, m_i \in \mathbb{Z}$. Then we have 
$$\varphi g(P) = \tilde{g}\varphi(P).$$
That is,  the following  diagram commutes
    \begin{figure}[h]
    \centering
	   \includegraphics[scale=1.2]{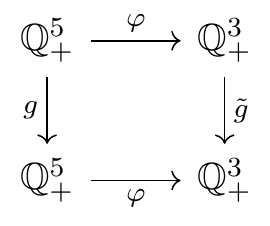}
    \end{figure}
\end{lem}

\begin{proof}
    It is easy to check that
	\begin{eqnarray}
        \mathcal{C}_0  \alpha = \mathcal{C}_0,~ \mathcal{C}_1  \alpha = \mathcal{C}_2,~ \mathcal{C}_2  \alpha = \mathcal{C}_0 \cdot \mathcal{C}_2 -\mathcal{C}_1,\label{alphaCi}\\
        \mathcal{C}_0  \beta = \mathcal{C}_1,~ \mathcal{C}_1  \beta = \mathcal{C}_2,~ \mathcal{C}_2  \beta = \mathcal{C}_0.\label{betaCi}
    \end{eqnarray}
	Then we have	
	\begin{eqnarray*}
		 	\varphi(\alpha(P)) &=& (\mathcal{C}_0(\alpha(P)), \mathcal{C}_1(\alpha(P)), \mathcal{C}_2(\alpha(P)))\\
					&=&  (\mathcal{C}_0(P), \mathcal{C}_2(P), \mathcal{C}_0(P)\mathcal{C}_2(P)-\mathcal{C}_1(P))\\
					&=& \tilde{\alpha}(\mathcal{C}_0(P),\mathcal{C}_1(P),\mathcal{C}_2(P))\\
		 			&=& \tilde{\alpha}(\varphi(P)),
	\end{eqnarray*}
and
	\begin{eqnarray*}
		 \varphi(\beta(P)) &=& (\mathcal{C}_0(\beta(P)), \mathcal{C}_1(\alpha(P)), \mathcal{C}_2(\beta(P)))\\
					&=&  (\mathcal{C}_1(P), \mathcal{C}_2(P), \mathcal{C}_0(P))\\
					&=& \tilde{\beta}(\mathcal{C}_0(P), \mathcal{C}_1(P),\mathcal{C}_2(P))\\
		 			&=& \tilde{\beta}(\varphi(P)).
	\end{eqnarray*}
 
Similarly we can prove that $\varphi(\alpha^{-1}(P)) = \tilde{\alpha}^{-1}(\varphi(P)),  \varphi(\beta^{-1}(P)) = \tilde{\beta}^{-1}(\varphi(P))$. Hence, we have
    $$ \varphi(g(P)) = \varphi(\alpha^{n_1}\beta^{m_1} \cdots \alpha^{n_t}\beta^{m_t}(P)) = \tilde{\alpha}^{n_1}\tilde{\beta}^{m_1} \cdots \tilde{\alpha}^{n_t}\tilde{\beta}^{m_t}(\varphi(P)) = \tilde{g}(\varphi(P)).$$

\end{proof}

It is clear that $\widetilde{T}\tilde{g} = \widetilde{T}$ for any $\tilde{g} \in \widetilde{G}$. Combining this identity with the identity in \text{Proposition \ref{Tg=T}}, \text{Lemma \ref{TM2}} and \text{Lemma \ref{varphig}}, we have a commutative diagram showed in Figure \ref{commute}.\\

\begin{figure}[h]
\centering
	\includegraphics[scale=0.9]{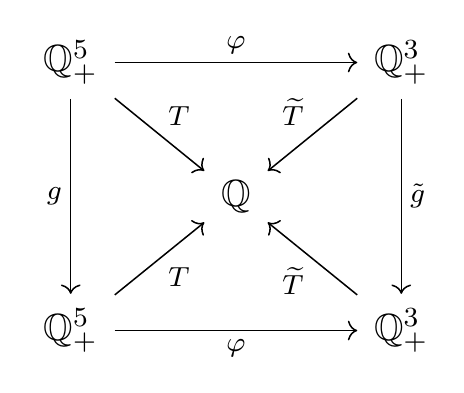}
\centering
\captionof{figure}{The commutative diagram of $ T, \widetilde{T}, g, \tilde{g}$ and $\varphi$.}\label{commute}
\end{figure}


    Now, we can prove that the set of the positive integral solutions of equation $\widetilde{T} = 0$ is the orbit of the initial solution $(3,4,4)$ under the group $\widetilde{G}$. We begin with a lemma.

\begin{lem}\label{perm} 
	If the triple $(x,y,z)$ belongs to the orbit $\widetilde{G}(3,4,4)$, so do the six permutations of $(x,y,z)$.
\end{lem}

\begin{proof}
	Let the triple $(x,y,z) = \tilde{g}(3,4,4)$ for some $\tilde{g} \in \widetilde{G}$. Since
	$(y,z,x) = \tilde{\beta}(x,y,z) $ and $(z,x,y)=\tilde{\beta}^{2}(x,y,z) $, we know $(y,z,x), (z,x,y) \in \widetilde{G}(3,4,4)$.
	
	Now we prove $(x,z,y) \in \widetilde{G}(3,4,4)$. Suppose $\tilde{g} = \gamma_n \cdots \gamma_1$, where $\gamma_i \in \{\tilde{\alpha}^{\pm1},  \tilde{\beta}^{\pm1}\}$. We denote a permutation $\sigma$ satisfying $\sigma(x,y,z) := (x,z,y)$. It is easy to check that $\sigma \gamma_i = \gamma_i^{-1}\sigma$. Hence 
	\begin{eqnarray*}
	    (x,z,y) &=& \sigma(x,y,z)\\
	            &=& \sigma\gamma_n \cdots \gamma_1(3,4,4)\\
	            &=& \gamma_n^{-1} \cdots \gamma_1^{-1}\sigma(3,4,4)\\
	            &=& \gamma_n^{-1} \cdots \gamma_1^{-1}(3,4,4) \in \widetilde{G}(3,4,4). 
	\end{eqnarray*}
\end{proof}

With above lemmas we can prove following theorem. 
\begin{thm} \label{S(344)}
	The set of all the positive integer solutions of the Diophantine equation $\widetilde{T}(X,Y,Z) = 0$ which defined in equation \eqref{markov+2} is the orbit of the initial solution $(3,4,4)$ under the group $\widetilde{G}$.
    That is, we have
	\begin{equation}
	\widetilde{G}(3,4,4) = \{(x,y,z) \in \mathbb{N}^3 \mid \widetilde{T}(x,y,z) = 0\}.
\end{equation}	
\end{thm}

\begin{proof} 
	``$\subseteq$'': Obviously, the triple $(3,4,4)$ is a positive integer solution. We claim that if the triple $(x,y,z)$ is a positive integer solution of the equation $\widetilde{T}(X,Y,Z)=0$, so does the triple $\tilde{g}(3,4,4)$ for all $\tilde{g} \in \widetilde{G}$. It is true, since 
    \begin{eqnarray}
	   \tilde{\alpha}(x,y,z) &=&  (x,z,xz-y) = (x, z, \frac{x^2 + z^2 + 7}{y}),\label{tilde}\\
	   \tilde{\alpha}^{-1}(x,y,z) &=& (x,xy-z,y) = (x, \frac{x^2 + y^2 + 7}{z},y),\label{tilde-}
	\end{eqnarray}
    and $\tilde{\beta}^{\pm1}(x,y,z)$ are positive integer solutions. Therefore, for all $\tilde{g} \in \widetilde{G}$, the triple $\tilde{g}(3,4,4)$ is a positive integer solution.

	``$\supseteq$'': Using \text{Lemma \ref{perm}}, we only need to consider the situation when $x\leq y \leq z$, where $(x,y,z)$ is a positive integer solution of the equation $\widetilde{T}(X,Y,Z)=0$.
	
	When $x < y = z$, we have $x^2 - xy^2+2y^2+7=0$. Since $x$ is an integer,  $\Delta = (y^2-4)^2-44$ must be a square number. By $(y^2-4+\sqrt{\Delta})(y^2-4-\sqrt{\Delta}) = 22 \times 2$ , we get $y = 4$. So $(x,y,z)=(3,4,4)$. When $x = y < z$, by the similar discussion, we have $(x,y,z) = (4,4,13)$ and we can check $\tilde{\alpha}\tilde{\beta}^2(3,4,4)=(4,4,13)$. When $x = y = z$, there is no integral solution. If $x=1$ or $2$, there is no integral solution. We know $x \geq 3, y \geq 4,z \geq 4$. 
	
	Now we prove the relationship of the two sets by induction on the maximum $m$ of triple $(x,y,z)$. When $m=4$, we have $(x,y,z)=(3,4,4) \in \widetilde{G}(3,4,4)$. 
	
	Assume that $m>4$. Then by above discussion, we know $y \ne z$. We consider $\tilde{\alpha}^{-1}(x,y,z) = (x, xy-z,y)$. Let $z' = xy-z$, it is a positive integer. We introduce a polynomial $f(\lambda) = \lambda^2-xy\lambda+x^2+y^2+7$. Note that $z, z'$ are zeros of $f$. 
	It is easy to check that $f(y)  = (3-x)y^2+(x+y)(x-y)+7 \leq 0$. If $f(y) = 0$, then $y=z'$. So $(x, z', y) = (3, 4, 4)$ and $z=8$. Hence $\tilde{\alpha}^{-1}(x,y,z)=\tilde{\alpha}^{-1}(3,4,8)=(3,4,4)$. If $f(y)<0$, then we know $z'<y<z$, which means the maximum of $(x, z', y)$, namely $y$, is smaller than the maximum of $(x, y, z)$, namely $z$. By induction hypothesis, there is an action $\tilde{h} \in \widetilde{G}$, such that $\tilde{h}(3,4,4) = (x, z', y)$. Then we have $\tilde{\alpha}\tilde{h}(3,4,4) = (x,y, z)$.
\end{proof}

\section{A description of the set $S(3,4,4)$ by action of a subgroup of $G$}\label{sec:pell}

    The last theorem urges us to consider the inverse image set $\varphi^{-1}(3,4,4)$. Additionally, from the observation of \text{Proposition \ref{mod}}, we also need to consider the set of vectors with the constant $\epsilon$ as the common multiple of each components. Therefore, we should determine the set $S(3,4,4)$, where
    \begin{equation} \label{Sci}
        S(C_0, C_1, C_2) := \varphi^{-1}(C_0, C_1, C_2) \cap  \{P=(a,b,c,d,\epsilon) \in \mathbb{N}^5 \mid P\equiv 0 (\textrm{mod~} \epsilon)\}.
    \end{equation}
    In the end of this section we will prove that the set $S(3,4,4)$ is the orbit of the initial solution $(\epsilon,\epsilon,\epsilon,\epsilon,\epsilon)$ under the subgroup $\langle\beta^3\rangle$, that is 
    $$S(3,4,4) = \langle\beta^3\rangle(\epsilon,\epsilon,\epsilon,\epsilon,\epsilon).$$
    
    We first begin with a lemma which gives a description of any vector in $\mathbb{Q}_+^5$. It can help us to determine the set $S(3,4,4)$.

\begin{lem} \label{repof5D}
	For any vector $P = (a,b,c,d,e) \in \mathbb{Q}_+^5$. Let a triple $(C_0, C_1, C_2) = \varphi(P)$, $m = C_0C_1C_2 - C_1^2 - C_2^2$ and $t = m - C_0^2 + 2$, where $\varphi$ defined in \eqref{varphi}. Then we have following three relations:
	
	(i) $$ c^2 - tcd + d^2 + C_0ce + C_0de + e^2 = 0$$
	
    (ii) $$P = \frac{1}{m}(c,d,e)\left(\begin{smallmatrix}
						{C_2} & {C_0C_2-C_1} & m  & 0 & 0\\
						{C_0C_1-C_2} & {C_1} & 0 & m & 0 \\
						{C_1} & {C_2} & 0 & 0 & m
					\end{smallmatrix}\right)$$
					
	(iii) $\forall n \in \mathbb{Z}$, 
	            $$\beta^{3n}(P) = P(L_{\text{sgn}(n)}(P))^{|n|},$$
    where $L_0(P) := \text{diag}(1,1,1,1,1)$,
    $$L_1(P) := \frac{1}{m}\left(\begin{smallmatrix}
						0&0&0&0&0\\
						0&0&0&0&0\\
						tC_2+C_0C_1-C_2 & t(C_0C_2-C_1)+C_1 & tm & m & 0\\
						-C_2 & C_1 -C_0C_2 & -m & 0 & 0\\
						C_1 - C_0C_2 & C_2 + C_0C_1 - C_0^2C_2 & -C_0m & 0 & m
					\end{smallmatrix}\right),$$
    $$L_{-1}(P) := \frac{1}{m}\left(\begin{smallmatrix}
						0&0&0&0&0\\
						0&0&0&0&0\\
						C_2 -C_0C_1 & -C_1 & 0 & -m & 0\\
						t(C_0C_1-C_2)+C_2 & tC_1+C_0C_2-C_1 & m & tm & 0\\
						C_1 + C_0C_2 - C_0^2C_1 & C_2 - C_0C_1 & 0 & -C_0m & m
					\end{smallmatrix}\right).$$
     
\end{lem}

   \begin{proof}

    (i) We first claim that the value $m$ is positive. It is true. Since by \text{Lemma \ref{TM2}} we have
	\begin{eqnarray*}
        m &=& \widetilde{T}(C_0,C_1,C_2) + C_0^2 + 7\\
		  &=& \widetilde{T}(\varphi(P)) + C_0^2 + 7\\
		  &=& T(P) + C_0^2 + 7\\
          &=& \frac{ab(c^2+d^2+e^2) + (a^2+b^2+cd)(c+d)e}{abcd} - 9 + \left( \frac{a^2 + b^2 + cd}{ab} \right)^2 +7\\
		  &=& \frac{ab(c^2+d^2+e^2) + (a^2+b^2+cd)(c+d)e}{abcd} + \frac{a^4+b^4+c^2d^2+2a^2cd+2b^2cd}{(ab)^2}\\
		  &>& 0.
	\end{eqnarray*}

	By above computation, we have
	\begin{eqnarray*}
		tcd &=& (m-C_0^2 + 2)cd\\
		    &=& (T(P)+9)cd\\
			&=& c^2 + d^2 + e^2 + \frac{a^2+b^2+cd}{ab}(c+d)e\\
			&=& c^2 + d^2 + e^2 + C_0ce + C_0de.
	\end{eqnarray*}
	Hence $c,d,e$ satisfy following equation
	$$c^2 -(m-C_0^2 + 2)cd+d^2+C_0ce+C_0de + e^2 = 0.$$


    (ii) Now we prove the second equation. We can check
	$$a = \frac{C_2c + (C_0C_1-C_2)d+C_1e}{m}, \quad b = \frac{(C_0C_2-C_1)c+C_1d+C_2e}{m}.$$
	
	
	Hence 
	$$P = (a,b,c,d,e) = \frac{1}{m}(c,d,e)\left(\begin{smallmatrix}
						{C_2} & {C_0C_2-C_1} & m  & 0 & 0\\
						{C_0C_1-C_2} & {C_1} & 0 & m & 0 \\
						{C_1} & {C_2} & 0 & 0 & m
					\end{smallmatrix}\right).$$

		
	
	(iii) If $n=0$, it is true.

    If $n>0$, since the value of $C_i$ will change after action by $\beta$, so we denote a matrix as follow 
	\begin{equation*}
	    M(P) := \left( \begin{smallmatrix}
					0 & 0 & 1 & 1 & 0\\
					1 & 0 & -\mathcal{C}_0(P) & 0 & 0\\    
					0 & 1 & \mathcal{C}_2(P) & 0 & 0\\
					0 & 0 & 0 & 0 & 0\\
					0 & 0 & 0 & 0 & 1
					\end{smallmatrix}\right).
	\end{equation*}
	
	We can  easily check that $ {(ac+be)}/{d} = a - \mathcal{C}_0(P)b + \mathcal{C}_2(P)c$. 
	Hence we have
	\begin{eqnarray*} 
	\beta(P) &=& \beta(a,b,c,d,e)\\
	        &=& (b,c,{(ac+be)}/{d}, a, e)\\
	        &=& (b, c, a -\mathcal{C}_0(P)b + \mathcal{C}_2(P)c, a, e)\\
	        &=& PM(P).
	\end{eqnarray*}
	
    By \text{the equation \eqref{betaCi}} we have
	\begin{eqnarray*}
	\beta^3(P) &=& \beta(\beta^2(P))\\
                &=& \beta^2(P)  M(\beta^2(P))\\
                &=& \beta(P)  M(\beta(P))  M(\beta^2(P))\\
                &=& P  M(P)  M(\beta(P))  M(\beta^2(P))\\
                &=& P  \left( \begin{smallmatrix}
					0 & 0 & 1 & 1 & 0\\
					1 & 0 & -\mathcal{C}_0(P) & 0 & 0\\    
					0 & 1 & \mathcal{C}_2(P) & 0 & 0\\
					0 & 0 & 0 & 0 & 0\\
					0 & 0 & 0 & 0 & 1
					\end{smallmatrix}\right)  \left( \begin{smallmatrix}
					0 & 0 & 1 & 1 & 0\\
					1 & 0 & -\mathcal{C}_0(\beta(P)) & 0 & 0\\    
					0 & 1 & \mathcal{C}_2(\beta(P)) & 0 & 0\\
					0 & 0 & 0 & 0 & 0\\
					0 & 0 & 0 & 0 & 1
					\end{smallmatrix}\right) \left( \begin{smallmatrix}
					0 & 0 & 1 & 1 & 0\\
					1 & 0 & -\mathcal{C}_0(\beta^2(P)) & 0 & 0\\    
					0 & 1 & \mathcal{C}_2(\beta^2(P)) & 0 & 0\\
					0 & 0 & 0 & 0 & 0\\
					0 & 0 & 0 & 0 & 1
					\end{smallmatrix}\right)\\
				&=& P  \left( \begin{smallmatrix}
					0 & 0 & 1 & 1 & 0\\
					1 & 0 & -C_0 & 0 & 0\\    
					0 & 1 & C_2 & 0 & 0\\
					0 & 0 & 0 & 0 & 0\\
					0 & 0 & 0 & 0 & 1
					\end{smallmatrix}\right)  \left( \begin{smallmatrix}
					0 & 0 & 1 & 1 & 0\\
					1 & 0 & -C_1 & 0 & 0\\    
					0 & 1 & C_0 & 0 & 0\\
					0 & 0 & 0 & 0 & 0\\
					0 & 0 & 0 & 0 & 1
					\end{smallmatrix}\right)  \left( \begin{smallmatrix}
					0 & 0 & 1 & 1 & 0\\
					1 & 0 & -C_2 & 0 & 0\\    
					0 & 1 & C_1 & 0 & 0\\
					0 & 0 & 0 & 0 & 0\\
					0 & 0 & 0 & 0 & 1
					\end{smallmatrix}\right)\\
				&=& P \left( \begin{smallmatrix}
					1 & {C}_0 & C_0C_1-C_2& 0 & 0\\
					-C_0 & 1-C_0^2 & C_1 + C_0(C_2-C_0C_1) & 0 & 0\\    
					 C_2 & C_0C_2 -C_1 & 1-C_1^2 + C_2(C_0C_1-C_2) & 1 & 0\\
					0 & 0 & 0 & 0 & 0\\
					0 & 0 & 0 & 0 & 1
					\end{smallmatrix}\right).
	\end{eqnarray*}
	
	Then by (ii), we have
	\begin{eqnarray*}
	\beta^3(P) &=& \frac{P}{m}\left( \begin{smallmatrix}
	    			0 & 0 & 0\\
					0 & 0 & 0\\
					1 & 0 & 0 \\ 
					0 & 1 & 0 \\
					0 & 0 & 1
			\end{smallmatrix}\right) \left(\begin{smallmatrix}
						{C_2} & {C_0C_2-C_1} & m  & 0 & 0\\
						{C_0C_1-C_2} & {C_1} & 0 & m & 0 \\
						{C_1} & {C_2} & 0 & 0 & m
					\end{smallmatrix}\right) \left( \begin{smallmatrix}
					1 & {C}_0 & C_0C_1-C_2& 0 & 0\\
					-C_0 & 1-C_0^2 & C_1 + C_0(C_2-C_0C_1) & 0 & 0\\    
					 C_2 & C_0C_2 -C_1 & 1-C_1^2 + C_2(C_0C_1-C_2) & 1 & 0\\
					0 & 0 & 0 & 0 & 0\\
					0 & 0 & 0 & 0 & 1
					\end{smallmatrix}\right)\\
			&=& \frac{P}{m}\left(\begin{smallmatrix}
						0&0&0&0&0\\
						0&0&0&0&0\\
						tC_2+C_0C_1-C_2 & t(C_0C_2-C_1)+C_1 & tm & m & 0\\
						-C_2 & C_1 -C_0C_2 & -m & 0 & 0\\
						C_1 - C_0C_2 & C_2 + C_0C_1 - C_0^2C_2 & -C_0m & 0 & m
					\end{smallmatrix}\right)\\
			&=& {P}L_1(P).
	\end{eqnarray*}

	On the other hand, we know that $\varphi(\beta^3(P)) = \tilde{\beta}^3(\varphi(P)) = (C_0, C_1, C_2)$, hence we have $L_1(\beta^3(P)) = L_1(P)$. Then we have
	$$\beta^{3n}(P) = \beta^3(\beta^{3(n-1)}(P)) = \beta^{3(n-1)}(P)L_1(\beta^{3(n-1)}(P)) = \beta^{3(n-1)}(P)L_1(P) = PL_1(P)^n. $$

    If $n < 0$, the proof is similar to the above case.
\end{proof}
    
    When we fix $(C_0, C_1, C_2) = (3,4,4)$ and $e = \epsilon$, the relation (i) of the last lemma inspire us to denote the Diophantine equation 
	\begin{equation}\label{Hyperbola}
	    H(X,Y) := X^2 - 9XY +Y^2+3\epsilon X+ 3\epsilon Y + \epsilon^2 = 0.
	\end{equation}
	We introduce a liner transformation  
	\begin{equation}
	    \theta(x) = 7x - 3\epsilon.
	\end{equation}
    And let $U = \theta(X), V = \theta(Y)$, then the equation \eqref{Hyperbola} is reduced to the hyperbola
\begin{equation}\label{hyperbola}
	\hat{H}(U,V) :=  U^2 - 9UV + V^2 +112\epsilon^2 = 0.
\end{equation}

{
   We know that the set of all integer solutions of a hyperbola can be obtained by Matthew's method\cite{Matthews}.
   
\begin{thm}[Matthews \cite{Matthews}]\label{matthew} 
    Given a hyperbola  
    \begin{equation}\label{ABCE}
	        AU^2 +BUV +CV^2 = E,
    \end{equation}
where $A,B,C,E$ are integers, $A>0, E < 0$ and $D = B^2 -4AC > 0$ is squarefree. Let $(x, y)$ be the least positive solution of Pell's equation
    \begin{equation}\label{eq:pell}
	X^2 - DY^2  = 4.
    \end{equation}
    Then the set of all the integer solutions of \eqref{ABCE} can be written as 
\begin{equation*}
	\bigcup_{(u,v)}\left\{ (u, v) \left(\begin{smallmatrix} (x-By)/2 & Ay  \\ -Cy & (x+By)/2 \end{smallmatrix}\right)^n \mid n \in \mathbb{Z}\right\},
\end{equation*}
    where $(u,v)$ is a positive integral solution of \eqref{ABCE} and satisfied one of the following conditions:
	 
	\quad (i) $\sqrt{4A|E|/D} \leq v < \sqrt{A|E|(x+2)/D}$,
	 
	\quad (ii) $u = \frac{1}{2A}(1-\frac{By}{x-2})\sqrt{A|E|(x-2)}$, $v = \sqrt{A|E|(x+2)/D}$.\\
\end{thm}}


	Note that our aim is to find the positive integer solutions of equation \eqref{Hyperbola} with $\epsilon$ as the common multiple of each component. So, we denote the set 
	\begin{equation}\label{solHyper}
		H_0 := \{ (x,y) \in \mathbb{N}^2 \mid H(x,y) = 0, (x,y)\equiv 0 (\textrm{mod~} \epsilon) \}.
	\end{equation}
	Applying above theorem, we obtain the set $H_0$ a precise description.

\begin{lem}\label{H0} For the equation $H(X,Y)=0$ which denoted in \eqref{Hyperbola}, we can obtain the set $H_0$ of the positive integer solutions of this equation with $\epsilon$ as the common multiple of each component, defined in \eqref{solHyper}. That is,
	\begin{equation}
	H_0 = \{ (x, y) \mid (x, y, \epsilon) = (\epsilon, \epsilon, \epsilon) \left(\begin{smallmatrix}9 & 1&0 \\ -1 &0&0\\ -3 &0 &1\end{smallmatrix}\right)^n , n \in \mathbb{Z} \}.
	\end{equation}
\end{lem}

\begin{proof} 
    We use three steps to prove it. 
    
    \textbf{First step.} We determine the set $\hat{H}_0$, where \begin{equation*}
		\hat{H}_0 := \{ (x,y) \in \mathbb{N}^2 \mid \hat{H}(x,y) = 0, (x,y)\equiv 0 (\textrm{mod~} \epsilon) \}.
	\end{equation*}
    and $\hat{H}$ is the hyperbola defined in \eqref{hyperbola}, that is, $\hat{H}(U,V) =  U^2 - 9UV + V^2 +112\epsilon^2$.

    Applying \text{Theorem \ref{matthew}}, we know that the set of all the positive integer solutions of the hyperbola \eqref{hyperbola} is 
\begin{equation*}
	\bigcup_{(u,v)}\left\{ (u, v) \left(\begin{smallmatrix} 9 & 1  \\ -1 & 0 \end{smallmatrix}\right)^n \mid n \in \mathbb{Z}\right\},
\end{equation*}
    where $(u,v)$ is a positive integral solution of \eqref{hyperbola} that satisfies 
    \begin{equation}\label{cond}
        \frac{8}{\sqrt{11}}\epsilon \leq v < 4\epsilon \quad \text{or} \quad (u,v) = (32\epsilon, 4\epsilon).       
    \end{equation}
    
    Only when  $(u,v)$ equals $(32\epsilon, 4\epsilon)$ does it satisfy the condition \eqref{cond} and  $(u,v) \equiv 0 (\textrm{mod~} \epsilon)$. 
    And, it is easy to check that $(u,v) \in \hat{H}_0$ if and only if $(u, v) \left(\begin{smallmatrix} 9 & 1  \\ -1 & 0 \end{smallmatrix}\right)^n \in \hat{H}_0$ for any $n \in \mathbb{Z}$. Hence, we have a precise description of $\hat{H_0}$, 
	\begin{equation} \label{hatH0}
		\hat{H}_0 = \{  (32\epsilon, 4\epsilon) \left(\begin{smallmatrix} 9 & 1  \\ -1 & 0 \end{smallmatrix}\right)^n \mid n \in \mathbb{Z} \}.
	\end{equation}
	
	\textbf{Second step.} We use $\hat{H}_0$ to determine $H_0$. We claim that 
	\begin{equation}\label{h0}
		H_0 = \{ (\theta^{-1}(u), \theta^{-1}(v)) \mid (u,v) \in \hat{H}_0\}.
	\end{equation}	
	
	``$\subseteq$'': Trivial.
	
	``$\supseteq$'': 
	Let $(u,v) = (32\epsilon, 4\epsilon) \left(\begin{smallmatrix} 9 & 1  \\ -1 & 0 \end{smallmatrix}\right)^{n}$  for some $n \in \mathbb{Z}$, we should prove $(\theta^{-1}(u),\theta^{-1}(v)) \in H_0$. Without lose of generality, we assume $n \geq 0$, and we prove it by induction on exponent $n$. When $n=0$, clearly $(\theta^{-1}(32\epsilon) , \theta^{-1}(4\epsilon)) = (5\epsilon, \epsilon) \in H_0$.
	
	Assume it is true for $n$, that is, suppose $(\theta^{-1}(u),\theta^{-1}(v)) \in H_0$.
	Since 
	$$(32\epsilon, 4\epsilon)\left(\begin{smallmatrix} 9 & 1  \\ -1 & 0 \end{smallmatrix}\right)^{n+1} = (u,v)\left(\begin{smallmatrix} 9 & 1  \\ -1 & 0 \end{smallmatrix}\right) = (9u-v, u),$$
	we should prove  $(\theta^{-1}(9u-v), \theta^{-1}(u)) \in H_0$.
	On one hand, we note that $w \equiv 4\epsilon (\textrm{mod~} 7\epsilon)$ if and only if  $\theta^{-1}(w) \equiv 0 (\textrm{mod~} \epsilon)$. 
	Hence  by $(\theta^{-1}(u),\theta^{-1}(v)) \equiv 0 (\textrm{mod~} \epsilon)$, we know  $(u,v) \equiv 4\epsilon (\textrm{mod~} 7\epsilon)$.
	Then we can check $ (9u-v, u) \equiv 4\epsilon (\textrm{mod~} 7\epsilon)$,  this implies  $(\theta^{-1}(9u-v), \theta^{-1}(u)) \equiv 0 (\textrm{mod~} \epsilon)$.
	On the other hand, it is easy to check that $(\theta^{-1}(9u-v), \theta^{-1}(u))$ is integral solution of the equation $H(X,Y)=0$. Hence we prove that $(\theta^{-1}(9u-v), \theta^{-1}(u)) \in H_0$.
	
	\textbf{Third step.} we simplify the set $H_0$. Let $(u,v) \in \hat{H}_0$. Since
	\begin{eqnarray*}
(\theta^{-1}(u),\theta^{-1}(v),\epsilon) &=& (u,v,\epsilon)\left(\begin{smallmatrix} 7 & 0&0 \\  0& 7&0\\ -3 & -3&1\end{smallmatrix}\right)^{-1}\\
	&=& (32\epsilon, 4\epsilon, \epsilon) \left(\begin{smallmatrix} 9 & 1 &0 \\ -1 & 0&0\\ 0&0&1 \end{smallmatrix}\right)^n   \left(\begin{smallmatrix} 7 & 0&0 \\  0& 7&0\\ -3 & -3&1\end{smallmatrix}\right)^{-1} \hspace{1.5cm}\text{(By \eqref{hatH0})}\\
	&=& (4\epsilon, 4\epsilon, \epsilon) \left(\begin{smallmatrix} 9 & 1 &0 \\ -1 & 0&0\\ 0&0&1 \end{smallmatrix}\right)^{n+1}   \left(\begin{smallmatrix} 7 & 0&0 \\  0& 7&0\\ -3 & -3&1\end{smallmatrix}\right)^{-1}\\
	&=& (\epsilon, \epsilon, \epsilon)\left(\begin{smallmatrix} 7 & 0&0 \\  0& 7&0\\ -3 & -3&1\end{smallmatrix}\right) \left(\begin{smallmatrix} 9 & 1 &0 \\ -1 & 0&0\\ 0&0&1 \end{smallmatrix}\right)^{n+1}\left(\begin{smallmatrix} 7 & 0&0 \\  0& 7&0\\ -3 & -3&1\end{smallmatrix}\right)^{-1} \\
	&=&  (\epsilon,\epsilon,\epsilon)\left(\begin{smallmatrix}9 & 1&0 \\ -1 &0&0\\ -3 &0 &1\end{smallmatrix}\right)^{n+1}
\end{eqnarray*}
for some $n \in \mathbb{Z}$, we can simplify the set $H_0$ in the equation \eqref{h0} as 
	\begin{equation*}
	H_0 = \{ (x, y)  \mid (x, y, \epsilon) = (\epsilon, \epsilon, \epsilon) \left(\begin{smallmatrix}9 & 1&0 \\ -1 &0&0\\ -3 &0 &1\end{smallmatrix}\right)^n , n \in \mathbb{Z} \}.
	\end{equation*}
\end{proof}

    Recall that $S(3,4,4) = \varphi^{-1}(3,4,4) \cap  \{P=(a,b,c,d,\epsilon) \in \mathbb{N}^5 \mid P\equiv 0 (\textrm{mod~} \epsilon)\}$, now we give this set a precise description.
\begin{thm} \label{beta3} The set $S(3,4,4)$ defined in \eqref{Sci} is the orbit of the vector $(\epsilon,\epsilon,\epsilon,\epsilon,\epsilon)$ under the subgroup $\langle\beta^3\rangle$ defined in \eqref{eq: alpha beta}, that is 
    $$S(3,4,4) = \langle\beta^3\rangle(\epsilon,\epsilon,\epsilon,\epsilon,\epsilon).$$
\end{thm} 
\begin{proof}  
	``$\supseteq$'': By \text{Lemma \ref{varphig}}, for any $n \in \mathbb{Z}$ we have
	\begin{eqnarray*}
	\varphi(\beta^{3n}(\epsilon,\epsilon,\epsilon,\epsilon,\epsilon)) &=& \tilde{\beta}^{3n}(\varphi(\epsilon,\epsilon,\epsilon,\epsilon,\epsilon)) \\
	&=&  \tilde{\beta}^{3n}(3,4,4) \\
	&=& (3,4,4).
	\end{eqnarray*}
	And by \text{Proposition \ref{mod}}, we have $\beta^{3n}(\epsilon,\epsilon,\epsilon,\epsilon,\epsilon) \in S(3,4,4)$.\\
	
	``$\subseteq$'': 
    Let $P = (a,b,c,d,\epsilon) \in S(3,4,4)$. Then there exists $n \in \mathbb{Z}$, such that
\begin{eqnarray*}
	    P &=& (c, d, \epsilon)\left(                 \begin{smallmatrix}
			1/4 & 2/4 & 1  & 0 & 0\\
            2/4 &1/4  & 0 & 1 & 0 \\
			1/4 & 1/4 & 0 & 0 & 1
   		    \end{smallmatrix}\right) \hspace{3.5cm}\text{(Lemma \ref{repof5D}.(ii))}\\
         &=&  (\epsilon, \epsilon, \epsilon)        \left(\begin{smallmatrix}
			9 & 1&0 \\
		    -1 &0&0\\
			-3 &0 &1
            \end{smallmatrix}\right)^n
			\left( \begin{smallmatrix}
        		1/4 & 2/4 & 1  & 0 & 0\\
			    2/4 &1/4  & 0 & 1 & 0 \\
				1/4 & 1/4 & 0 & 0 & 1
			\end{smallmatrix}\right) \hspace{1.8cm}\text{(Lemma \ref{H0})}\\
\end{eqnarray*}

If $n=0$, we have $P = (\epsilon, \epsilon, \epsilon, \epsilon, \epsilon)$.

If $n>0$, then 
\begin{eqnarray*}     
     P &=&   (\epsilon, \epsilon, \epsilon, \epsilon, \epsilon) \left(\begin{smallmatrix}
			0 & 0&0 &0 & 0\\
			0 & 0&0 &0 & 0\\
			0 & 0 &9 & 1&0 \\
			0 & 0 & -1 &0&0\\
			0 & 0 &  -3 &0 &1\end{smallmatrix}\right)^n \left( \begin{smallmatrix}
			0 & 0&0 &0 & 0\\
			0 & 0&0 &0 & 0\\
			1/4 & 2/4 & 1 & 0&0 \\ 
			2/4 & 1/4 & 0& 1 & 0\\
			1/4 & 1/4 & 0& 0& 1
    \end{smallmatrix}\right)^n\\
    		& = & (\epsilon,\epsilon,\epsilon,\epsilon,\epsilon) 	\left( \begin{smallmatrix}
		0 & 0 & 0 & 0 & 0\\
		0 & 0 & 0 & 0 & 0\\    
		11/4 & 19/4 & 9 & 1 & 0\\
		-1/4 & -2/4 & -1 & 0 & 0\\
		-2/4 & -5/4 & -3 & 0 & 1
	\end{smallmatrix}\right)^n\\
	& = & (\epsilon,\epsilon,\epsilon,\epsilon,\epsilon)(L_1(\epsilon,\epsilon,\epsilon,\epsilon,\epsilon))^n\\
	& = & \beta^{3n}(\epsilon,\epsilon,\epsilon,\epsilon,\epsilon). \hspace{5cm}\text{(Lemma \ref{repof5D}.(iii))}
\end{eqnarray*}

If $n<0$, then 
\begin{eqnarray*}     
     P &=&   (\epsilon, \epsilon, \epsilon, \epsilon, \epsilon) \left(\begin{smallmatrix}
			0 & 0&0 &0 & 0\\
			0 & 0&0 &0 & 0\\
			0 & 0 & 0 & 1 & 0 \\
			0 & 0 & 1 & 9 & 0\\
			0 & 0 & 0 & -3 &1\end{smallmatrix}\right)^{-n} \left( \begin{smallmatrix}
			0 & 0&0 &0 & 0\\
			0 & 0&0 &0 & 0\\
			1/4 & 2/4 & 1 & 0&0 \\ 
			2/4 & 1/4 & 0& 1 & 0\\
			1/4 & 1/4 & 0& 0& 1
    \end{smallmatrix}\right)^{-n}\\
    		& = & (\epsilon,\epsilon,\epsilon,\epsilon,\epsilon) 	\left( \begin{smallmatrix}
		0 & 0 & 0 & 0 & 0\\
		0 & 0 & 0 & 0 & 0\\    
		-2/4 & 1/4 & 0 & -1 & 0\\
		-19/4 & 11/4 & 1 & 9 & 0\\
		-5/4 & -2/4 & 0 & -3 & 1
	\end{smallmatrix}\right)^{-n}\\
	& = & (\epsilon,\epsilon,\epsilon,\epsilon,\epsilon)(L_{-1}(\epsilon,\epsilon,\epsilon,\epsilon,\epsilon))^{-n}\\
	& = & \beta^{3n}(\epsilon,\epsilon,\epsilon,\epsilon,\epsilon). \hspace{5cm}\text{(Lemma \ref{repof5D}.(iii))}
\end{eqnarray*}    
\end{proof}

\section{The main result and its proof}\label{sec:bao}

	Depending on the above two sections, now we can answer the Lampe's question: Fix the constant $\epsilon \in \mathbb{N}$, which solutions $(a,b,c,d,\epsilon) \in \mathbb{N}^5$ of the Diophantine equation $T(a,b,c,d,\epsilon)=0$ can be obtained from the initial solution $(\epsilon,\epsilon,\epsilon,\epsilon,\epsilon)$ by a sequence of actions of the group $G$? The answer is quite simple: 
	$$\textrm{if and only if~} \varphi(a,b,c,d,\epsilon) \in \mathbb{N}^3 \textrm{~and~} (a,b,c,d,\epsilon) \equiv 0 (\textrm{mod~}\epsilon).$$
	
	We begin with a lemma.

\begin{lem}\label{gP}
	Given a vector $P =  (a,b,c,d,\epsilon) \in \mathbb{N}^5$. If the triple $ \varphi(P) \in \mathbb{N}^3$ and $P \equiv 0 (\textrm{mod~}\epsilon)$, then for any action $g \in G$, we have $ \varphi(g(P)) \in \mathbb{N}^3$ and $ g(P) \equiv 0 (\textrm{mod~}\epsilon) $.
\end{lem}
\begin{proof}
	Denote the triple $(C_0, C_1, C_2) = \varphi(P)$. By \text{Lemma \ref{varphig}}, it is easy to check that $\varphi(\beta^{\pm 1}(P)) = \tilde{\beta}^{\pm 1}(\varphi(P)) \in \mathbb{N}^3$. By the equations \eqref{tilde} and \eqref{tilde-} we know $\varphi(\alpha^{\pm 1}(P)) = \tilde{\alpha}^{\pm 1}(\varphi(P)) \in \mathbb{N}^3$.
	
	Since $P \equiv 0 (\textrm{mod~}\epsilon)$ and 
	\begin{eqnarray*}
		\alpha(P) &=& (b, \frac{b^2+cd}{a}, c,d,\epsilon) = (b, C_0b-a,c,d,\epsilon),\\
		\alpha^{-1}(P) &=& (\frac{a^2+cd}{b},a,c,d,\epsilon) = (C_0a-b,a,c,d,\epsilon),\\
		\beta(P) &=& (b,c,\frac{ac+b\epsilon}{d},a,\epsilon) = (b,c,a-C_0b+C_2c,a,\epsilon),\\
		\beta^{-1}(P) &=&(d, a, b, \frac{a\epsilon+bd}{c}, \epsilon) = (d, a, b, -C_0a+b+C_1d, \epsilon),
	\end{eqnarray*}
	we know  $\alpha^{\pm 1}(P) \equiv \beta^{\pm 1}(P) \equiv 0 (\textrm{mod~}\epsilon) $. Hence for any $g \in G$, we have $ \varphi(g(P)) \in \mathbb{N}^3, g(P) \equiv 0 (\textrm{mod~}\epsilon)$.
\end{proof}

    Now we prove the main result.
\begin{thm}\label{MainThm} Fix the constant $\epsilon \in \mathbb{N}$. The positive integer solutions $(a,b,c,d,\epsilon) \in \mathbb{N}^5$ of the Diophantine equation $T(a,b,c,d,\epsilon)=0$ defined in \eqref{T=9} can be obtained from the initial solution $(\epsilon,\epsilon,\epsilon,\epsilon,\epsilon)$ by a sequence of actions of the group $G:= \langle\alpha,\beta\rangle$ defined in \eqref{eq: alpha beta}, if and only if, $\varphi(a,b,c,d,\epsilon) \in \mathbb{N}^3 \textrm{~and~} (a,b,c,d,\epsilon) \equiv 0 (\textrm{mod~}\epsilon)$, where the map $\varphi$ defined in \eqref{varphi}. That is, we have 
    $$G(\epsilon,\epsilon,\epsilon,\epsilon,\epsilon) = \{ P=(a,b,c,d,\epsilon) \in \mathbb{N}^5 \mid T(P) = 0, \varphi(P) \in \mathbb{N}^3, P \equiv 0
	(\textrm{mod~}\epsilon) \}.$$
\end{thm}

\begin{proof}
    ``$\subseteq$'': For any action $g \in G$, denote a vector $P = g(\epsilon,\epsilon,\epsilon,\epsilon,\epsilon)$. First, by  \text{Proposition \ref{Tg=T}}, we have
	$$T(P) = T(g(\epsilon,\epsilon,\epsilon,\epsilon,\epsilon)) = T(\epsilon,\epsilon,\epsilon,\epsilon,\epsilon) = 0.$$	

	Second, by \text{Lemma \ref{varphig}} and \text{Theorem \ref{S(344)}}, we have
    $$\varphi(P)  = \varphi (g(\epsilon,\epsilon,\epsilon,\epsilon,\epsilon)) =  \tilde{g}(\varphi(\epsilon,\epsilon,\epsilon,\epsilon,\epsilon)) = \tilde{g}(3,4,4) \in \mathbb{N}^3.$$
	
	Last, by  \text{Proposition \ref{mod}}, we have
	$$P  \in \mathbb{N}^5 \quad \textrm{and} \quad P \equiv 0 (\textrm{mod~} \epsilon).$$

	Hence $P$ belongs to the right set.\\

	``$\supseteq$'': Let $P=(a,b,c,d,\epsilon) \in \mathbb{N}^5$ belongs to the right set. Denote $(C_0, C_1, C_2) = \varphi(P)$. By \text{Lemma \ref{TM2}},we have 
	\begin{equation*}
		\widetilde{T}(C_0, C_1, C_2) =  \widetilde{T}(\varphi(P))= T(P) = 0.
	\end{equation*}
So, by \text{Theorem \ref{S(344)}}, there exists an action $\tilde{g} \in \widetilde{G}$ such that $\tilde{g}(3,4,4)=(C_0, C_1, C_2)$. Suppose $\tilde{g} = \tilde{\alpha}^{n_1}\tilde{\beta}^{m_1}\cdots\tilde{\alpha}^{n_t}\tilde{\beta}^{m_t}$ for some $n_i, m_i \in \mathbb{Z}$. Let $g = {\alpha}^{n_1}{\beta}^{m_1}\cdots {\alpha}^{n_t}{\beta}^{m_t}$. On the one hand, by \text{Lemma \ref{varphig}}, we have
	\begin{eqnarray*}
	    \varphi(g^{-1}(P)) &=& \varphi({\beta}^{-m_t}{\alpha}^{-n_t}\cdots{\beta}^{-m_1}{\alpha}^{-n_1}(P))\\
                &=& \tilde{\beta}^{-m_t}\tilde{\alpha}^{-n_t}\cdots\tilde{\beta}^{-m_1}\tilde{\alpha}^{-n_1}(\varphi(P))\\
                &=& \tilde{\beta}^{-m_t}\tilde{\alpha}^{-n_t}\cdots\tilde{\beta}^{-m_1}\tilde{\alpha}^{-n_1}(C_0, C_1, C_2)\\
                &=& \tilde{\beta}^{-m_t}\tilde{\alpha}^{-n_t}\cdots\tilde{\beta}^{-m_1}\tilde{\alpha}^{-n_1}\tilde{g}(3,4,4)\\
                &=& (3,4,4).
	\end{eqnarray*}
    On the other hand, by  \text{Lemma \ref{gP}}, we know that $ \varphi(g^{-1}(P)) \in \mathbb{N}^3$ and $ g^{-1}(P) \equiv 0 (\textrm{mod~}\epsilon)$. Therefore  $g^{-1}(P) \in S(3,4,4)$.
    
	Then by \text{Theorem \ref{beta3}}, there exist an integer $n$, such that $\beta^{3n}(\epsilon,\epsilon,\epsilon,\epsilon,\epsilon) = g^{-1}(P)$. Hence
	$$P = g\beta^{3n}(\epsilon,\epsilon,\epsilon,\epsilon,\epsilon).$$
\end{proof}
\vspace{4mm}

{\bf Acknowledgements:}\;  {\em This project is supported by the National Natural Science Foundation of China  (No.12071422 and No.12131015).}

\bibliographystyle{plain}

\begin{thebibliography}{1}

\bibitem{FZ1}
Sergey Fomin and Andrei Zelevinsky.
\newblock Cluster algebras. {I}. {F}oundations.
\newblock {\em J. Amer. Math. Soc.}, 15(2):497--529, 2002.

\bibitem{Hone}
Andrew N.~W. Hone.
\newblock Laurent polynomials and superintegrable maps.
\newblock {\em SIGMA Symmetry Integrability Geom. Methods Appl.}, 3:Paper 022,
  18, 2007.

\bibitem{lampe}
Philipp Lampe.
\newblock Diophantine equations via cluster transformations.
\newblock {\em J. Algebra}, 462:320--337, 2016.

\bibitem{Markov}
A.~Markoff.
\newblock Sur les formes quadratiques binaires ind\'{e}finies.
\newblock {\em Math. Ann.}, 17(3):379--399, 1880.

\bibitem{Matthews}
Keith~R. Matthews, John~P. Robertson, and Anitha Srinivasan.
\newblock On fundamental solutions of binary quadratic form equations.
\newblock {\em Acta Arith.}, 169(3):291--299, 2015.

\bibitem{PengZhang}
X.~Peng and J.~Zhang.
\newblock Cluster algebras and markoff numbers.
\newblock {\em CaMUS}, 3:19--26, 01 2012.

\end{thebibliography}

\end{document}